\documentclass[12pt]{iopart}
\usepackage{graphicx}
\usepackage[utf8]{inputenc}
\usepackage{setstack,iopams}
\usepackage{comment}

\usepackage{xcolor}
\usepackage[urlcolor=blue]{hyperref}      
\hypersetup{
  colorlinks = true,                    
  citecolor = {purple},
  linkcolor = {blue},
}

\newcommand{\erm}{\mathrm{e}}
\newcommand{\irm}{\mathrm{i}}
\renewcommand{\epsilon}{\varepsilon}
\newcommand{\qty}[1]{\left( #1 \right)}

\begin{document}
\title[Distributed time delays]{Mapping dynamical systems with distributed time delays 
to sets of ordinary differential equations}
\author{Daniel Henrik Nevermann$^{1,*}$, Claudius Gros$^1$}
\address{$^1$Institute for Theoretical Physics, 
Goethe University Frankfurt, Germany}
\ead{nevermann@itp.uni-frankfurt.de}

\begin{abstract}
Real-world dynamical systems with retardation effects 
are described in general not by a single, precisely
defined time delay, but by a range of delay times.
An exact mapping onto a set of $N+1$ ordinary 
differential equations exists when the respective
delay distribution is given in terms of a gamma
distribution with discrete exponents. The number 
of auxiliary variables one needs to introduce, $N$, 
is  inversely proportional to the variance of the 
delay distribution. The case of a single delay is
therefore recovered when $N\to\infty$. Using this 
approach, denoted here the `kernel series framework', 
we examine systematically how the bifurcation phase
diagram of the Mackey-Glass system changes under the
influence of distributed delays. We find that local
properties, f.i.\ the locus of a Hopf bifurcation, 
are robust against the introduction of broadened 
memory kernels. Period-doubling transitions and the 
onset of chaos, which involve non-local properties 
of the flow, are found in contrast to be more 
sensitive to distributed delays. In general,
the observed effects are found to scale as $1/N$.
Furthermore, we consider time-delayed systems 
exhibiting chaotic diffusion, which is present
in particular for sinusoidal flows. We find that
chaotic diffusion is substantially more pronounced
for distributed delays. Our results indicate
in consequence that modeling approaches of real-world 
processes should take the effects of distributed 
delay times into account.
\end{abstract}
\vfill
{\footnotesize$^*$ Corresponding author.}
\maketitle

\section{Introduction}

On a fundamental level, all known laws of
nature are strictly markovian. An example 
is the Schrödinger equation,
$\hbar\dot\psi(\mathbf{x},t)=H \psi(\mathbf{x},t)$,
which describes a complex-valued dynamical system 
for which the time evolution is fully determined by 
the current state of the system, the wavefunction 
$\psi(\mathbf{x},t)$. There is no memory. In
contrast, delayed effects often emerge in the
context of macroscopic processes. The study of 
effective models containing time delays is hence an important
subject \cite{lakshmanan2011dynamics}.

The vast majority of studies dedicated to time delayed 
systems assume that the delay $T$ is accurately defined
\cite{richard2003time,otto2019nonlinear},
with $T$ being either constant or functionally 
dependent on time, e.g.\ containing 
a periodic modulation \cite{muller2018laminar}.
Strictly speaking, time delays are however never precisely
defined. Consider a dynamical variable $x(t)$
whose time evolution is influenced by past states, 
say by $x(t-T)$. For this to be possible, the 
past trajectory must be stored, either explicitly
or implicitly, via suitable physical, chemical, or 
biophysical processes. Memory formation takes 
however time, which implies that only a washed-out 
version of the past will be available. Mathematically,
the system is hence described not by a fixed time
delay, but by distributed time delays.

The broadening of the  memory kernel may be disregarded
when it is small, which is achievable in particular for 
technical systems \cite{erneux2017introduction}.
For natural systems \cite{cooke1982discrete},
distributed time delays are however prominent. This
is well known for epidemic spreading, which can be modeled
on a realistic level only with distributed delays
\cite{beretta2001global,mccluskey2010complete}.
The underlying reason is that biological processes,
like incubation and recovery times, are intrinsically variable.
A corresponding observation holds for the interaction between
tumors and the immune systems \cite{sardar2021impact},
or for the delayed responses intrinsic to
predator-prey models \cite{xu2012bifurcations}.

It is well-known that dynamical systems with 
distributed time delays can be mapped onto sets 
of ordinary differential equations (ODEs) using
the so-called ``linear chain trick'' 
(see, e.g., \cite{gedeon1999upper, hurtado2020procedure,
mocek2005approximation, hurtado2019generalizations}).
The resulting ODE system takes a particular 
simple form when the distribution of the delays 
is described by a gamma distribution with integer
exponents. This method, which we denote here the 
``kernel series framework'', allows for the systematic 
study of delay differential equations (DDEs) 
with distributed delays and the influence
of broadened memory kernels. It is closely 
related to previous work in the
context of integro-differential equations \cite{valani2022lorenz}, non-reducible
distributed delays \cite{macdonald1987stability} and threshold delay systems
\cite{smith1993threshold}. 

For the case of a DDE with discrete time delays, 
the kernel series framework provides an
approximation scheme via a set of $N\!+\!1$ ODEs.
This approximation is systematic in the sense
that the original DDE is recovered in the limit
$N\!\to\!\infty$, with the differences scaling 
as $1/N$. The kernel series framework allows
hence to study delay systems with publicly 
available scientific software libraries, which
contain solvers typically only for sets of
ordinary differential equations, but not for 
delay systems.

The method is reviewed in section\,\ref{sec:theory},
where we show that a finite order $N$ of the kernel 
series corresponds to a delay distribution with a relative width 
$1/\sqrt{N}$. Subsequently, in section\,\ref{sec:results}, 
we study analytically and numerically first the
delay-induced Hopf bifurcation and then the
properties of a prototypical delay dynamical system, 
the Mackey-Glass system. Comparing the original DDE 
with the corresponding kernel series framework, we 
find distinct differences between the local and 
the global regime. For a good description of the Hopf 
bifurcation, which depends on local properties of 
the flow, only a modest number of auxiliary 
variables is needed, which implies that local 
bifurcations are robust with respect to
washed-out memory kernels.

Period-doubling transitions and delay-induced chaos
are phenomena that depend on non-local properties of
the flow. In this regime, the global regime, in part
substantial differences between DDEs with a single time 
delay and systems with distributed delays are observed.
Even for larger $N$, of the order of a few hundred,
both the respective bifurcation points and the 
topology of the resulting attractors may differ
noticeably. This result has important consequences
for applications. Globally stabilized dynamical 
states observed for delay systems in the limiting 
case of precisely defined time delays may not
be present in the corresponding real-world systems
characterized by distributed time delays.

In section \ref{sect:generic}, we discuss the 
extension of the kernel series framework to 
general time-delay distributions. Afterwards, 
in section \ref{sec:chaotic_diffusion}, we
consider an application to a system exhibiting
diffusive chaotic dynamics, finding that chaotic 
diffusion is substantially enhanced for 
distributed delays.

\section{Kernel series framework}\label{sec:theory}

The dynamics of a dynamical system with 
time delays depends on its history
$h(t)$, as defined by
\begin{equation}
\label{eq:general_dde}
\dot{x}(t) = F\qty{t, x(t), h(t)}\,,
\end{equation}
where $x(t)$ is the primary dynamical variable
and $F$ a given function.
For a DDE with a single time delay $T$ one has
$h(t)=x(t-T)$. In general, the memory $h(t)$ 
is given by a superposition of past states 
$x(t-\tau)$,
\begin{equation}
h(t) = \int_0^\infty x(t-\tau) K(\tau)\,\mathrm{d}\tau,
\qquad\quad
\int_0^\infty K(\tau)\,\mathrm{d}\tau=1\,,
\label{eq:timeDelayKernel}
\end{equation}
where $K(\tau)\ge0$ is the time delay kernel.
A single time delay of size $T$ is present for
$K(\tau)=\delta(\tau-T)$.

\begin{figure}[t]
\centering
\includegraphics[]{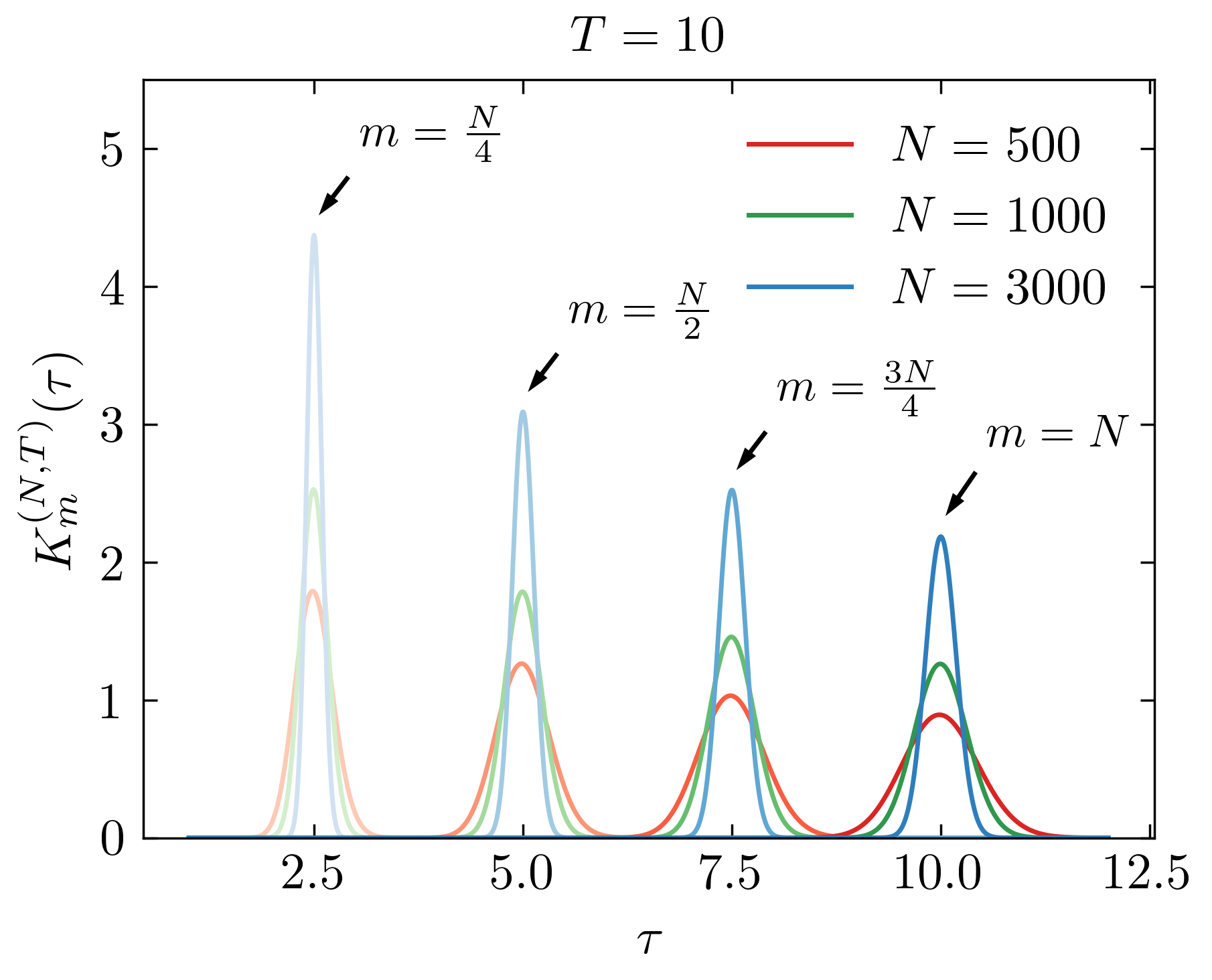}
\caption{{\bf Time delay kernels.} 
The probability densities of $K^{(N,T)}_{m}(\tau)$,
as defined by (\ref{eq:kernel_gamma}),
for $T=10$, together with
$N=500$ (red lines),
$N=1000$ (green lines) and
$N=3000$ (blue lines).
The four peaks observable for each $N$ correspond
respectively to $m/N=1.0,\,0.75,\,0.5,\,0.25$ 
(progressively shaded lines correspond
to distinct $m/N$). 
The width of the peaks scale as $1/\sqrt{m}$, see
(\ref{eq_mean_variance}).
The kernels converge to the Dirac delta distributions
$\delta(\tau-mT/N)$ in the limit $N\to\infty$.
}
\label{fig:kernel_convergence}
\end{figure}


We concentrate mostly on time delay kernels 
based on the probability density of the gamma 
distribution. Later on, in section\,\ref{sect:generic},
the general case will be treated. For a given 
number $N>0$, the order of the framework, 
we define $N$ kernels
\begin{equation}
K^{(N, T)}_{m}(\tau) = \frac{N^m}{(m-1)!\,T^m}\,
\left[\tau^{m-1}\erm^{-N \tau / T}\right],
\qquad\quad m = 1, \dots, N\,,
\label{eq:kernel_gamma}
\label{eq:kernel_def}
\end{equation}
which correspond to a set of $N$ normalized 
gamma distributions. Mean and variance are
\begin{equation}
\mu=\frac{m T}{N}, 
\qquad\quad
\sigma^2 = \frac{m T^2}{N^2},
\qquad\quad
\frac{\sigma}{\mu} = \frac{1}{\sqrt{m}}\,,
\label{eq_mean_variance}
\end{equation}
which implies that
\begin{equation}
\lim_{N\to\infty} K_N^{(N,T)}(t) = \delta(t-T)\,.
\label{eq:N_infty_delta}
\end{equation}
An illustration of the kernel series is given in
figure \,\ref{fig:kernel_convergence}. One observes
that the densities converge to a $\delta$-function
for $N\to\infty$ when the ratio $m/N$ is kept fixed,
in agreement with (\ref{eq_mean_variance}).
For a given memory $h(t)$, distributed by a gamma-shaped
delay kernel $K(\tau) \propto \tau^x \mathrm{e}^y$,
the corresponding mean time delay $T$ and order $N$
are easily retrieved from the inverse relations
\[
  T = \frac{x + 1}{y}, \qquad\quad N = x + 1\,.
\]

\subsection{Equivalence to sets of ordinary differential equations}

There are many representations of the Dirac $\delta$-functions commonly used.
What makes (\ref{eq:kernel_def}) especially interesting is that the individual
kernels can be evaluated by a simple additional ordinary differential equation,
which is the essence of the linear chain trick \cite{hurtado2020procedure,
mocek2005approximation, hurtado2019generalizations, cassidy2021distributed}. For
this purpose auxiliary variables $x_m$ are introduced via
\begin{equation}
x_m(t) = \int_0^\infty x(t-\tau) K_m^{(N,T)}(\tau)\,\mathrm{d}\tau, 
\qquad\quad
m = 1, \dots, N\,.
\label{eq_def_x_m}
\end{equation}
At face value, it may seem that one needs to store the full history of $x(t)$
for the evaluation of the convolutions defining the $x_m$. This is however not
the case. The reason is that the $x_m(t)$ form a closed set of recursive
differential equations \cite{gedeon1999upper, hurtado2019generalizations,
cassidy2021distributed},
\begin{equation}
\dot{x}_m(t) = \frac{x_{m-1}(t) - x_m(t)}{T_N},
\qquad\quad
x_0(t)=x(t), \qquad\quad T_N = \frac{T}{N}\,,
\label{eq_dot_x_m}
\end{equation}
which holds for $m=1,\,\dots,\,N$.
For a derivation one uses
\begin{eqnarray}
\label{eq:dot_m_derivation}
\dot{x}_m(t) &=&
\int_0^\infty K_m^{(N,T)}(\tau) \frac{\mathrm{d}}{\mathrm{d}t} x_0(t-\tau)\,\mathrm{d}\tau
\\
&=& -\int_0^\infty K_m^{(N,T)}(\tau) \frac{\mathrm{d}}{\mathrm{d}\tau} x_0(t-\tau)\,\mathrm{d}\tau\,
\nonumber
\end{eqnarray}
with an integration by parts of the last expression 
leading directly to (\ref{eq_dot_x_m}) (see e.g.\
\cite{hurtado2019generalizations} for a more rigorous proof).
Including $x_0(t)=x(t)$, one has then a dynamical 
system defined by $N\!+\!1$ variables
$x_m(t)$ ($m=0, \dots, N$).

\subsection{Kernel series dynamics}

The core object of interest is the 
$(N\!+\!1)$-dimensional dynamical system 
\begin{equation}
\label{eq:kernel_approx}
\cases{
\dot{x}_0(t) = F\qty{t, x_0(t), x_N(t)},
\qquad\quad x_0(t) = x(t),
	\\
\dot{x}_m(t) = \frac{x_{m-1}(t) - x_m(t)}{T_N},
\qquad\quad m = 1,\dots,N\,,
      } 
\end{equation}
where $T_N=T/N$. It can be viewed from two distinct
perspectives. First, that the
set of $N\!+\!1$ equations defined by
(\ref{eq:kernel_approx}) constitutes an exact mapping
of a DDE with distributed time delays to a set of
ordinary differential equations. This mapping holds 
when the respective time delay kernel is given
by a gamma distribution with an integer exponent.
We will show in section\,\ref{sect:generic} that
all delay distributions can be mapped to sets of ODEs 
when the number of auxiliary variables is 
correspondingly increased. 
For the most general case, a diverging number of 
memories $x_m(t)$ may however be necessary. 
Secondly, one can interpret (\ref{eq:kernel_approx}) 
as an approximation,
\begin{equation}
x(t-T) \approx x_N(t)\,,
\label{eq:x_x_N}
\end{equation}
to the delay system defined by (\ref{eq:general_dde})
when $h(t)=x(t-T)$. In this case the 
approximation becomes exact in the 
limit $N\!\to\!\infty$, see (\ref{eq:N_infty_delta}). 
For numerical investigations, one can use 
(\ref{eq:kernel_approx}) as a proxy for
the original DDE, with the advantage being
that standard ODE solvers can be used. In
the following we examine in detail what 
happens when the width of delay distributions
is increased or decreased, with distributions
of zero width corresponding to the case of
a single delay.

\subsection{Generalized state history / phase space collapse}

The trajectories of a delay system with 
a single delay $T$ are defined in the 
space of the state history \cite{wernecke2019chaos},
\begin{equation}
X(t) = \{x(t') \mid t' \in [t-T, t)\}\,,
\label{eq:state_history}
\end{equation}
which is formally infinite-dimensional.
Consistently, an initial function 
$\phi\colon[-T,0) \to \mathbb{R}$
is needed for the dynamics defined
by (\ref{eq:general_dde}) when
$h(t)=x(t-T)$. Interestingly, one can 
regard the set of dynamical variables 
$x_m(t)$ defining the kernel series
framework, (\ref{eq:kernel_approx}),
as a generalized state history. The dynamical 
history is however not sampled at specific points 
in time, as it would be the case for a discretized 
version of the standard state history. Instead, 
suitably weighted superpositions of the past 
trajectory are taken. 

The discretized version of the classical state 
history (\ref{eq:state_history}) can be used
to formulate a basic approximation to a DDE
in terms of discretized Euler updatings
\cite{wernecke2019chaos}. In contrast to the
ODE system (\ref{eq:kernel_approx}),
Euler updating schemes do however not incorporate a
systematic relation to distributed delay times.
The result, that the kernel series framework
is an exact representation of DDEs with 
gamma-distributed delay distributions with 
integer exponents, raises an interesting 
question. In general, delay systems
come with a formally infinite-dimensional 
state space, which does however collapse, as 
a corollary of above observation, for
specific delay distributions. An 
interesting point regards the condition
for this phase space collapse, namely
if it could occur also for other types
of delay distributions. We leave these 
questions for further investigations.

\section{Results}\label{sec:results}

We examined extensively the differences showing up
between the kernel series framework and the respective
original delay equation. We find that lower-order 
bifurcation transitions agree well already for 
modest $N$, which is however not the case for
chaotic attractors. In particular, we find that 
chaos tends to disappear when distributions of 
delays are considered, even when the width of 
delays is comparatively small.

For the presentation of the results we
distinguish between a local and a global 
regime, noting that the stability of a 
fixpoint is determined by the local
properties of the flow, with the properties
of chaotic attractors being determined by 
non-local, viz global properties of the
flow \cite{gros2015complex}. The same holds
for period-doubling transitions.

\begin{figure}[t]
\centering
\includegraphics[]{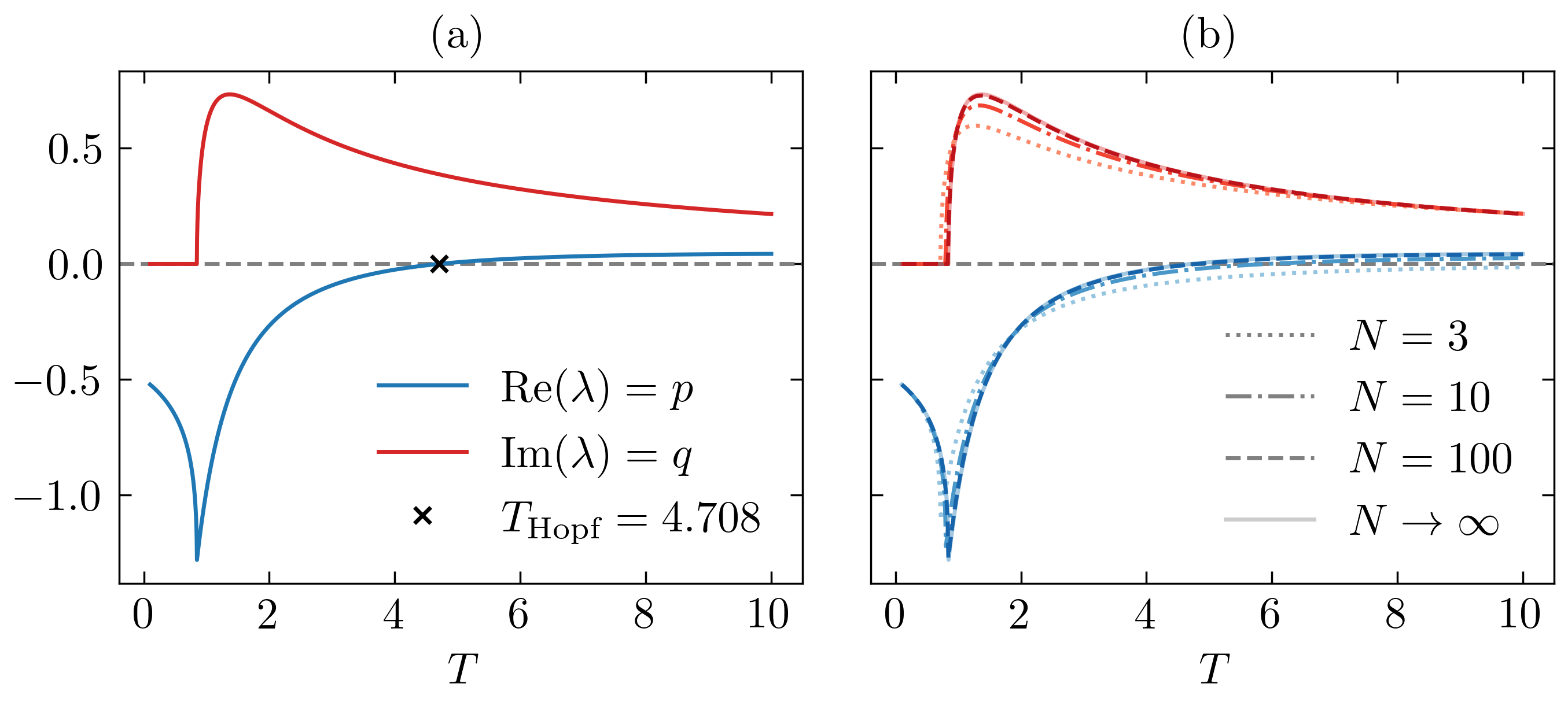}
\caption{{\bf Hopf bifurcation.}
\textbf{(a)} As a function of the time delay $T$, 
the real and imaginary part of the numerical solution
of (\ref{eq:characteristic_eq}) for
$\lambda=p + \irm q$, where the parameters of the linear delay system
(\ref{eq:linear_dde}) are $a = 0.4$ and $b = 0.1$.
A Hopf bifurcation occurs at $T_{\mathrm{Hopf}}$ when 
$\lambda$ crosses the imaginary axis, i.e.\ when 
$\mathrm{Re}(\lambda)=p=0$.
\textbf{(b)} The exponent $\lambda$ obtained from solving 
the analytic expression (\ref{eq:ode_characteristic_eq})
for the kernel series framework compared to the solution
shown in \textbf{(a)} obtained for the linear delay system
(\ref{eq:linear_dde}) which is found in the kernel series
framework by taking the limit $N\to\infty$. Shown are results for
$N=3,\,10,\,100$.
}
\label{fig:characteristic_eq}
\end{figure}

\subsection{Local regime: Hopf bifurcation}\label{sec:local_regime}

As a first example we study the generic delay system
\begin{equation}
\dot{x}(t) = -c - b x(t) - a x(t-T)\,,
\label{eq:linear_dde}
\end{equation}
which emerges when expanding (\ref{eq:general_dde})
with $h(t) = x(t-T)$ around a given fixpoint.
It is amenable to an analytic solution \cite{gros2015complex}.
To simplify discussions, we choose $c = 0$, which moves the
fixpoint of the system to $x^*= x(t) = x(t-T) = 0$.

The fixpoint $x^*$ is stable
for $b,\,a >0$ when the time delay $T$ is small.
The ansatz $x(t) = C\exp(\lambda t)$ leads to
\begin{equation}
\label{eq:characteristic_eq}
0= \lambda+b+a \erm^{-\lambda T}\,.
\end{equation}
A Hopf bifurcation occurs at a critical time 
delay $T_\mathrm{Hopf}$ when the real part of 
$\lambda=p+\irm q$ vanishes, viz when $\lambda$
crosses the imaginary axis. One finds
\begin{equation}
b = -a\cos(qT_{\mathrm{Hopf}}),
\quad
q = a \sin(qT_{\mathrm{Hopf}}),
\quad
T_{\mathrm{Hopf}} = \frac{\arccos{(-b/a)}}{\sqrt{a^2 - b^2}}.
\label{eq:hopf_T_c}
\end{equation}
A numerical solution of (\ref{eq:characteristic_eq}) 
is presented in figure\,\ref{fig:characteristic_eq}.
We use $a = 0.4$ and $b = 0.1$, which is consistent 
with the values of the linearized Mackey-Glass 
system, as discussed in section\,\ref{sec:global_regime}. The 
Hopf bifurcation point is then $T_{\mathrm{Hopf}} = 4.708$.

The kernel series framework (\ref{eq:kernel_approx}) for
the linear delay system (\ref{eq:linear_dde}) leads
to a linear system of $N+1$ coupled ordinary differential
equations. The fixpoint, which corresponds to $x_m=0$ for
$m=0,\dots,N$, has the Jacobian
\begin{equation}
  \mathbf{J} = \frac{1}{T_N}\left(\matrix{
    - b T_N &  &  & - a T_N \cr
    1 & - 1 &  &  \cr
      & \ddots & \ddots &  \cr
      &  & 1 & - 1
  }\right)\,, \qquad\quad T_N = \frac{T}{N}\, .
\label{eq:linear_Jacobian}
\end{equation}
The $N+1$ eigenvalues are retrieved from 
the roots of the characteristic polynomial
$P_N(\lambda)$, which can be extracted
analytically,
\begin{equation}
P_N(\lambda) = (b + \lambda) (1 + T_N \lambda)^N + a\,.
\label{eq:ode_characteristic_eq}
\end{equation}
As expected, the expression (\ref{eq:characteristic_eq})
for a single time delay $T$ is recovered
in the limit $N\to\infty$ when using
$\exp(x) = \lim_{k\to\infty} (1 + x/k)^k$ 
together with $T_N = T/N$. 

\begin{figure}[t]
\centering
\includegraphics[]{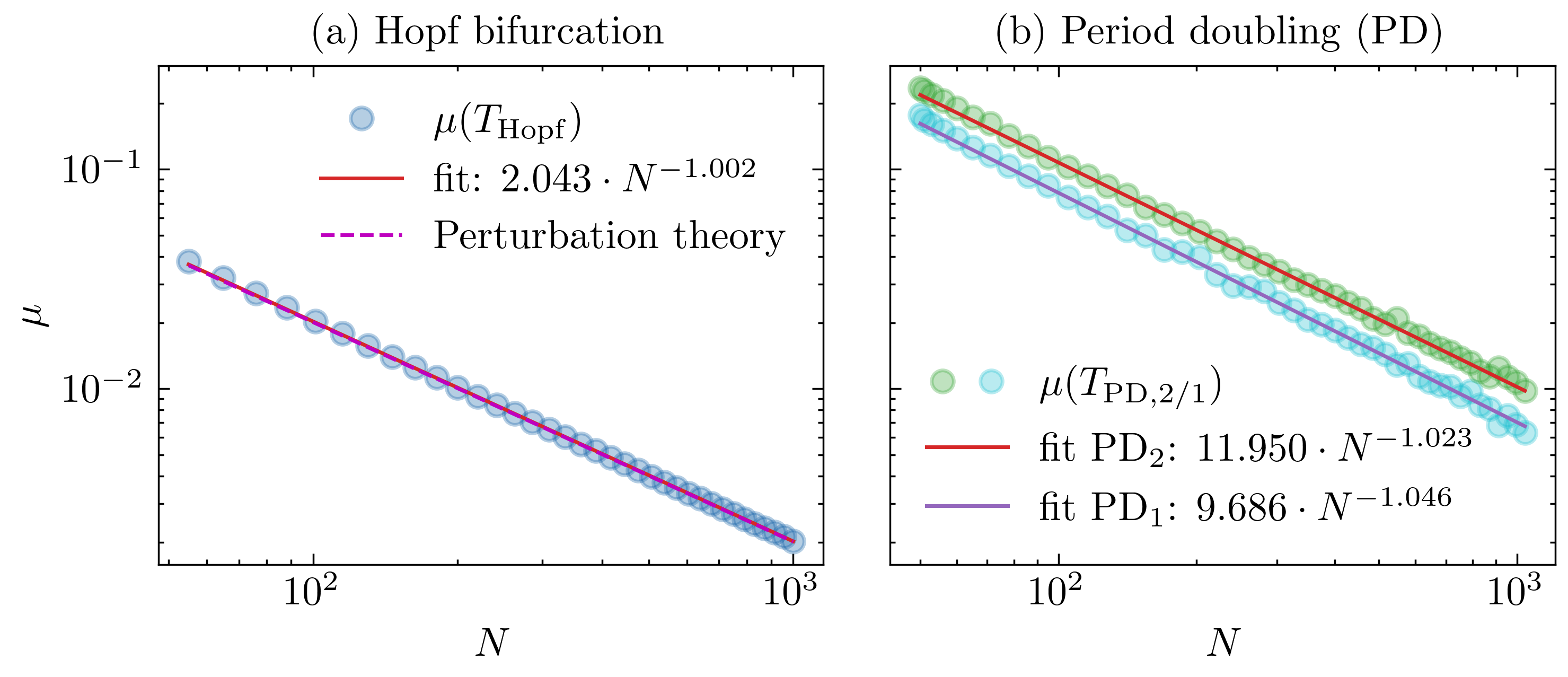}
\caption{
{\bf Relative differences.} Logarithmic representation 
of the relative differences between the Hopf bifurcation point
and the first and second period doubling bifurcation points
obtained within kernel series frameworks of order $N$,
(\ref{eq:kernel_approx}), and the value found by solving 
the delay system obtained correspondingly in the limit 
$N\!\to\!\infty$.
\textbf{(a)} For $a=0.4$ and $b=0.1$, the relative deviation
$\mu=\mu(T_\mathrm{Hopf})$ of the Hopf bifurcation 
point, as defined by (\ref{eq:mu_hopf}).
One observes excellent agreement between a linear fit 
(red line) and the analytic $1/N$ expansion (dashed magenta line), 
as given by (\ref{eq:appedix_T_Hopf}), which we derive in the 
\nameref{sec:appendix}.
\textbf{(b)} The relative deviation $\mu=\mu(T_\mathrm{PD})$ 
of the first (PD$_1$, blue circles) and second 
(PD$_2$, green circles) period doubling transitions, 
together with linear fits (purple/red lines). For a given
order $N$, the relative deviation is generally lowest for 
the Hopf bifurcation, and lower for the first period 
doubling than for the second.
}
\label{fig:hopf_rel_err}
\end{figure}

The dynamics of the system, and in particular its 
stability, is dictated by the maximal eigenvalue of
the Jacobian $\lambda_\mathrm{max}$.
In figure\,\ref{fig:characteristic_eq} we included
results for $\lambda_\mathrm{max}$
obtained by solving
(\ref{eq:ode_characteristic_eq}) numerically
for various orders $N$ of the kernel series framework.
The resulting maximal eigenvalues
converge quickly towards the ones found for 
the original system as $N$ increases,
as shown in figure\,\ref{fig:characteristic_eq}.

To further quantify the differences between
distributed and single time delays we evaluated 
the relative deviation of the $N$th 
order result for the Hopf bifurcation point,
$T^{(N)}_{\mathrm{Hopf}}$, from the 
value found for a single time delay (\ref{eq:hopf_T_c})
using
\begin{equation}
\label{eq:mu_hopf}
\mu(T_{\mathrm{Hopf}}) = 
\frac{\left| T_{\mathrm{Hopf}} - T^{(N)}_{\mathrm{Hopf}}\right|}
{T_{\mathrm{Hopf}}},
\qquad\quad
T_{\mathrm{Hopf}} =\lim_{N\to\infty} T^{(N)}_{\mathrm{Hopf}}
\end{equation}
The numerical results, as well as a power-law fit
to the data, are shown in figure\,\ref{fig:hopf_rel_err}\,(a).
We further compare to the analytical prediction 
attained via perturbation theory,
as presented in the \nameref{sec:appendix}. 
Asymptotically, the relative deviation 
$\mu(T_{\mathrm{Hopf}})$ scales as $1/N$, 
inversely with the order of the kernel 
series. Quantitatively an agreement of 5\% 
is achieved by $N\approx39$, and
1\% for $N\approx225$.

Here we also treat the notable case of exponentially
distributed delays that emerge for $N = 1$. The corresponding
characteristic polynomial only has solutions for negative
real part of the eigenvalue and thus no Hopf transition
is observed. This implies that for all $T$ the dynamics
are governed by the stable fixpoint.

\begin{figure}[t]
\centering
\includegraphics[]{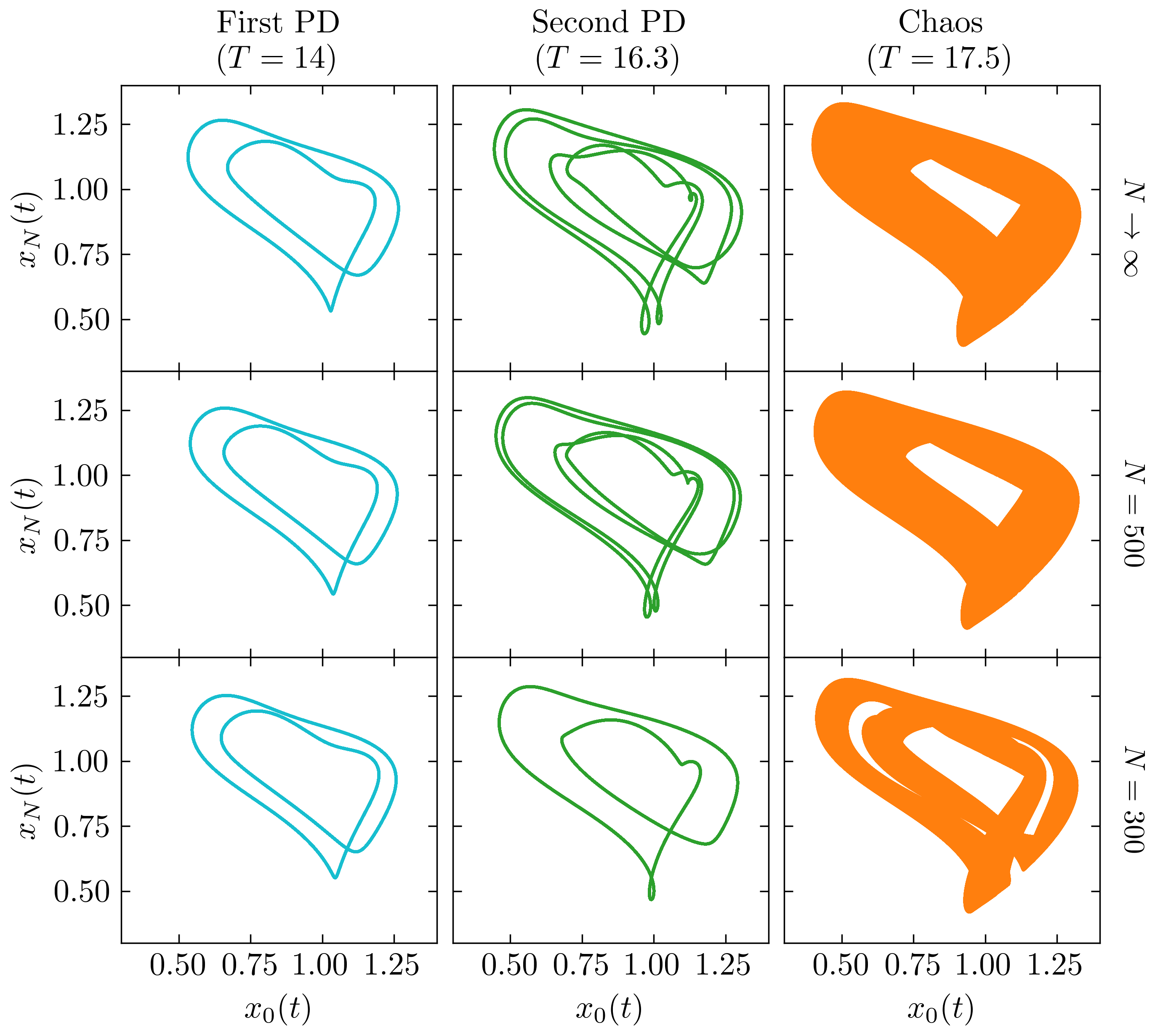}
\caption{{\bf Mackey-Glass system.}
The stroboscopic projection $(x_0(t),x_N(t))$ of attractors 
of the Mackey-Glass system (\ref{eq:mg}), where
$x_0(t)\to x(t)$ and $x_N\to x(t-T)$ for $N\to\infty$. Shown are numerical solutions of 
the original delay differential equation (top row), 
and the respective $N$th-order kernel framework, 
namely for $N = 500$ and $N=300$ (middle and bottom row). 
The average time delay $T$ is chosen 
such that the original system is beyond the first
period doubling (PD, left column), the second PD
(middle column), or in the chaotic regime (right column).
Substantial differences can be seen between the
system with a single time delay (top row) and
distributed delays when $N$ is not overly large
(bottom row).
}
\label{fig:mg_attractors}
\end{figure}

\subsection{Global regime}\label{sec:global_regime}

For the study of the influence of distributed
time delays in the global regime, we consider 
a prototypical time delay system with
chaotic attractors, the Mackey-Glass system,
\begin{equation}
\label{eq:mg}
\dot{x}(t) = \frac{\alpha x(t-T)}{1 + 
\qty{x(t-T)}^\gamma} - \beta x(t),
\quad\qquad \alpha, \beta, \gamma, T > 0\,,
\end{equation}
originally introduced to model the production 
of blood cells \cite{mackey1977oscillation}. 
In the following, we set the parameters to 
the standard values
$\alpha = 0.2$, $\beta = 0.1$, together
with $\gamma = 10$
\cite{lakshmanan2011dynamics,wernecke2019chaos}. 
This choice of parameters ensures that the 
trivial fixpoint $x_1^* = 0$ is unstable for all
time delays, while the fixpoint 
$x_2^* = (\alpha / \beta - 1)^{1/\gamma} = 1$ 
looses stability via a Hopf bifurcation when 
increasing the delay. In the following we 
examine to which extent the sequence of 
bifurcations occurring in (\ref{eq:mg})
changes when the relative width 
$\sigma/\mu=1/\sqrt{N}$,
see (\ref{eq_mean_variance}), of the delay 
distribution becomes positive. For an overview
of state-of-the-art studies of the Mackey-Glass 
see f.i.\ 
\cite{lakshmanan2011dynamics,wernecke2019chaos}. 


For small time delays, the dynamics are governed by the 
stable fixpoint $x_2^*$. Linearizing (\ref{eq:mg})
around $x_2^*$ and inserting an exponential ansatz yields
(\ref{eq:characteristic_eq}) with $b = 0.1$ and $a = 0.4$, 
the parameter values used in section\,\ref{sec:local_regime}. 
A supercritical Hopf bifurcation, destabilizing the
fixpoint in favor of a stable limit cycle occurs 
at $T_{\mathrm{Hopf}} \approx 4.708$, as 
determined by (\ref{eq:hopf_T_c}). When further 
increasing the time delay, the system undergoes a 
series of period doubling bifurcations 
(the first two at 
$T_{\mathrm{PD},1}\approx 13.39$ and
$T_{\mathrm{PD},2} \approx 15.95$)
and finally, for $T_{\mathrm{Chaos}} \approx 16.5$, 
a transition to a chaotic regime. We further note that 
the Mackey-Glass system is strictly dissipative, in
the sense that the divergence of the flow is always
negative. This holds also for the kernel series
framework, viz for distributed time delays, 
independently of the order $N$. 

In figure\,\ref{fig:mg_attractors} the stroboscopic
projection, plotting $x_N(t)\approx x(t-T)$ as a function
of $x_0(t)=x(t)$, is used to illustrate the 
topology of the attractors found for three 
different values of the average time delay, respectively
for $T_{\mathrm{PD},1} < T\!=\!14.0 <T_{\mathrm{PD},2}$
and $T_{\mathrm{PD},2} < T\!=\!16.3 <T_{\mathrm{PD},3}$,
as well as for $T=17.5$, which is a bit beyond 
the onset of chaos. Shown are the attractors obtained 
by solving (\ref{eq:mg}) directly, denoted as 
$N\!\to\!\infty$, together with the attractors 
obtained from the corresponding kernel 
series framework (\ref{eq:kernel_approx}),
both for $N=300$ and $500$.\footnote{The 
Mackey-Glass system is solved numerically using a 
software package introduced in \cite{jitcxde}, 
which employs a Runge-Kutta method, as described in
\cite{shampine2000solving}. The system of ordinary
differential equation generated by the kernel series framework
is solved using a standard fourth-order Runge-Kutta method.}

It is clear from figure\,\ref{fig:mg_attractors},
that the topology of attractors may change once
distributions of time delays with positive coefficient of variation 
are allowed. For example, for $T=16.3$, the resulting
limit cycle is doubled only once for $N=300$, but
twice for $N=500$ and above. Substantial differences
show up in addition for the chaotic attractors 
stabilizing for $T=17.5$ when the order $N$ is changed.
For $N=500$, the stroboscopic projection of the
chaotic attractor shown figure\,\ref{fig:mg_attractors},
seems to be similar, on first sight, to the 
$N\to\infty$ limit. However, we did not attempt 
to make a more precise comparison, f.i.\ in terms 
of the respective fractal dimensions, which we leave
for future studies.

\begin{figure}[t]
\centering
\includegraphics[]{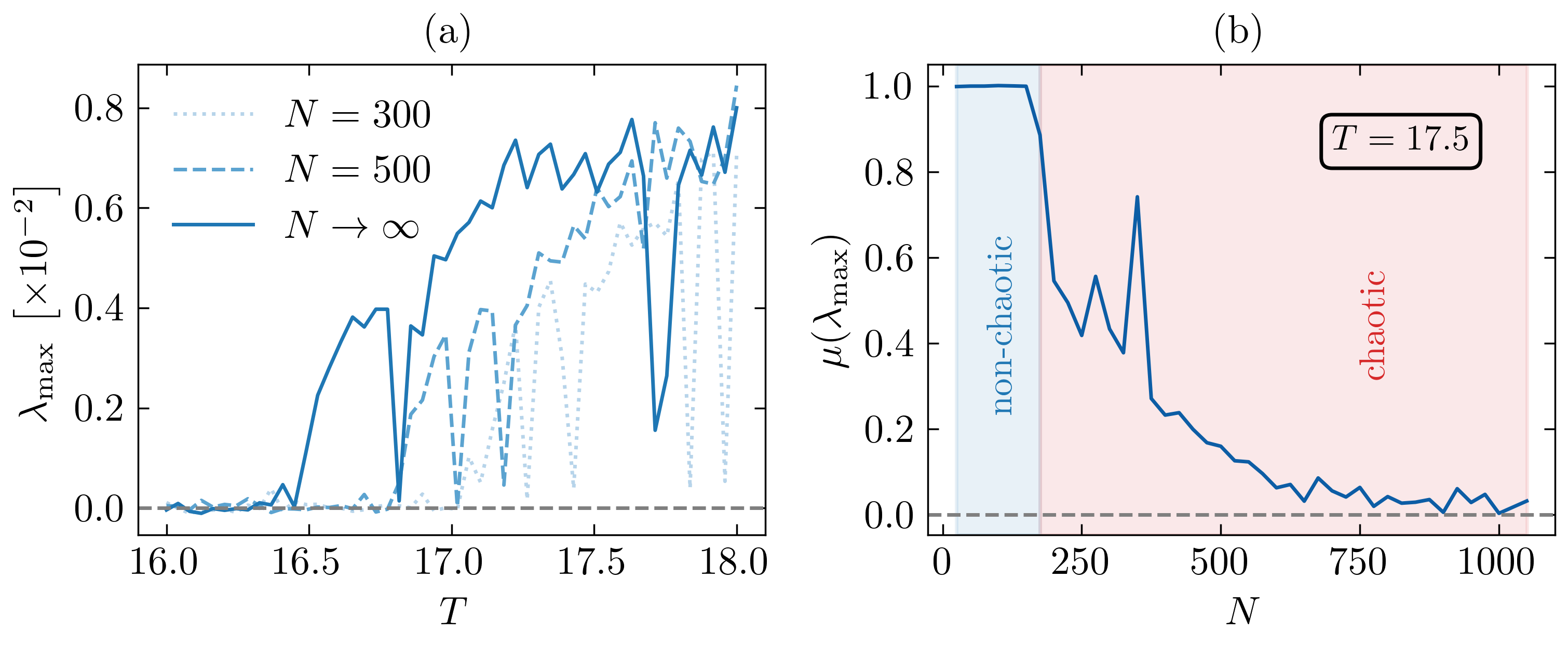}
\caption{
{\bf Lyapunov exponents.} 
The maximal Lyapunov exponent $\lambda_\mathrm{max}$ of the 
Mackey-Glass system, see (\ref{eq:mg}).
\textbf{(a)} As a function of the average time delay $T$, the
Lyapunov exponents obtained from the kernel series framework
(\ref{eq:kernel_approx}) for different orders,
$N = 300$ and $N=500$, with $N=\infty$ corresponding to
the case of a single delay.
\textbf{(b)} For $T = 17.5$, the relative deviation 
of the maximal Lyapunov exponent, obtained 
respectively for the given $N$ and for $N\to\infty$. 
For small $N$ (blue shading) the system has not 
yet entered the chaotic regime (red shading).
}
\label{fig:le}
\end{figure}

In order to quantify the differences between the kernel
series framework and the solution of (\ref{eq:mg}), we
consider in figure\,\ref{fig:hopf_rel_err}\,(b)
the locus of the bifurcation points of the first and 
second period doubling, $T_{\mathrm{PD}}$. Shown
are the relative differences $\mu$, as defined by 
(\ref{eq:mu_hopf}), between the results obtained 
numerically for finite $N$ and $N\to\infty$. One 
finds $1/N$ scaling, in analogy to the behavior observed 
for the primary Hopf bifurcation, as presented in
figure\,\ref{fig:hopf_rel_err}\,(a). Quantitatively,
$\mu=\mu(T_{\mathrm{PD}})$ drops below 5\% for the first
period doubling when $N>154$, and below 1\% for
$N>716$ kernels. For larger time delays, larger number
of kernels are required -- for the second period doubling
transition, $\mu$ drops below 5\% and 1\% respectively
for $N>212$ and $N>1020$.

Beyond a cascade of period doubling bifurcations,
the Mackey-Glass system enters a chaotic regime.
The most common measure for deterministic chaos 
in dynamical systems is the emergence of one or 
more positive Lyapunov exponents 
\cite{gros2015complex}. The dynamics is dominated 
by the maximal Lyapunov exponent $\lambda_{\mathrm{max}}$,
which quantifies the average spreading of two 
initially close-by trajectories. Systems with more
than one positive Lyapunov exponent are usually called
hyperchaotic \cite{rossler1979equation}. In the 
Mackey-Glass system, we observe a transition
to chaos for time delays $T > T_{\mathrm{Chaos}} \approx 16.5$. 
At this point, the maximal Lyapunov exponent
becomes positive, as shown in figure\,\ref{fig:le}\,(a).
Further on, for $T \gtrsim 27$, the Mackey-Glass system
shows hyperchaotic dynamics.

In figure\,\ref{fig:le}\,(a) the maximal Lyapunov exponent
of the Mackey-Glass system is compared to the maximal
Lyapunov exponents of the corresponding kernel series
frameworks. The point where the transition to chaos 
occurs decreases in general with increasing $N$.
For $N = 500$, we find chaotic behavior for 
$T \gtrsim 16.8$. In figure\,\ref{fig:le}\,(b), the relative
deviation of the maximal Lyapunov exponent of
the kernel series framework with respect to $N\to\infty$
is plotted over $N$. Systems with $N \lesssim 200$ 
have not yet transitioned to the chaotic regime at 
$T = 17.5$, which is associated with a kernel width of
$\sigma \gtrsim 1.3$. Thus, we note that chaos may break down
in the kernel series framework if the kernel width exceeds a
threshold value. In this sense, chaos is not robust
in the kernel series framework. The peak observed around
$N \approx 300$ is caused by a non-chaotic window within the
chaotic regime.

\subsection{Zero-One test for chaos}

The evolution of the cross-correlation of two initially 
close-by trajectories can be used to classify the
long-term dynamical behavior \cite{wernecke2017test}.
On defines with $C_{12}(t)$,
\begin{equation}
\label{eq:cross_correlation}
C_{12}(t) = 
\frac{\big\langle (x_1(t) - \bar{x}) (x_2(t) -\bar{x}) \big\rangle}{s^2} 
\equiv 1 - \frac{D_{12}(t)}{2 s^2},
\end{equation}
the cross-correlation of two trajectories
$x_1(t)$ and $x_2(t)$, where $\bar{x}$ is the
center of gravity of the attractor and $s$
the standard deviation. An average
$\langle\cdot\rangle$ over initial positions
$x_1(0)$ and $x_2(0)$ is performed, such that
the initial distance 
$\delta=\parallel\!x_1(0)-x_2(0)\!\parallel$ 
is kept constant, with $\parallel\!\cdot\!\parallel$ 
denoting the distance between the respective
initial functions.
The last term in (\ref{eq:cross_correlation})
connects $C_{12}(t)$ with the quadratic distance 
between the two trajectories,
$D_{12}(t) = \langle (x_1(t) - x_2(t))^2 \rangle$.

\begin{figure}[t]
\centering
\includegraphics[]{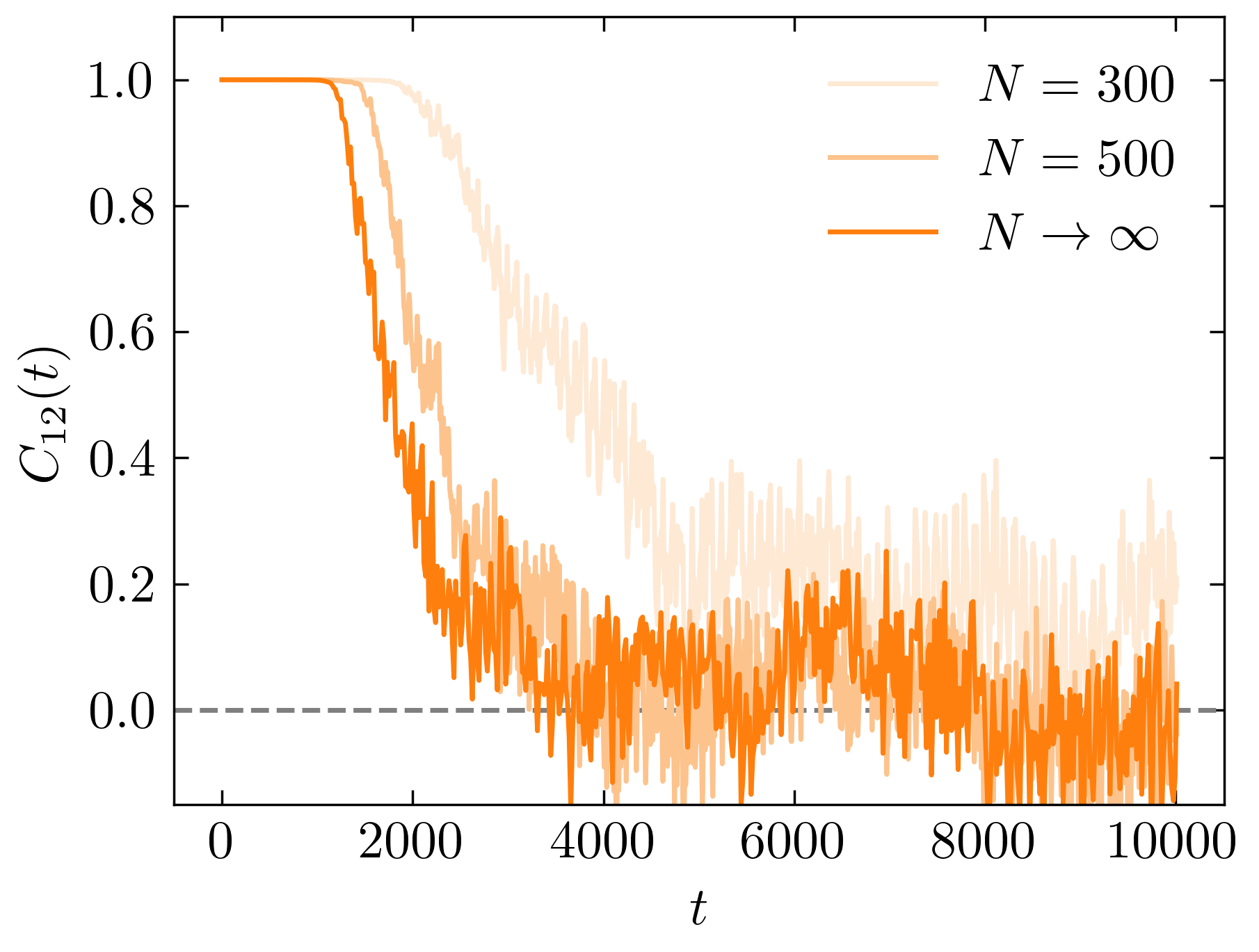}
\caption{\textbf{Cross-correlation.}
As a function of time $t$, the cross-correlation 
$C_{12}(t)$ of two initially close-by trajectories, 
as defined by (\ref{eq:cross_correlation}).
The plots are averages over 100 pairs of initial
functions with distances $\delta=10^{-6}$.
The average time delay is $T = 17.5$, for which
the Mackey-Glass system (\ref{eq:mg}) shows normal
turbulent chaos. Shown are the results of the
kernel series framework (\ref{eq:kernel_approx}) 
for $N = 300, 500$, together with the cross-correlation 
of the Mackey-Glass system with a single time delay,
denoted by $N \to \infty$.}
\label{fig:cross_correlation}
\end{figure}

The long-term behavior of chaotic and
non-chaotic dynamics differ qualitatively
with respect to $C_{12}$ and $D_{12}$, which 
can be used hence as a One-Zero test for chaos.
For the four basic types of attractors
one has \cite{wernecke2017test}:
\begin{itemize}
\item {\em Fixpoint:} A fixpoint attracting both
trajectories leads to $D_{12}\to0$ 
independently of the initial distance $\delta$.

\item {\em Limit Cycle:} On average, two 
trajectories ending up in the same limit cycle 
have a finite distance that scales with the
initial distance. For a limit cycle one has
hence $D_{12}\propto\delta$.

\item {\em Chaotic Attractor:} Independently
of the initial distance, trajectories become fully 
decorrelated on a chaotic attractor.
This implies that $D_{12}\propto s^2$ for any
$\delta>0$.

\item {\em Partially Predictable Chaos:} In this
case the initial exponential divergence is limited
by topological constraints, with the final chaotic
state being reached only via a subsequent diffusive 
process, which can be very slow. One has $D_{12}< s^2$ 
for an extended period.

\end{itemize}
All four states are found in the Mackey-Glass 
system \cite{wernecke2017test,wernecke2019chaos}.
In figure\,\ref{fig:cross_correlation} 
the cross-correlation is shown for 
an average delay $T = 17.5$. For large
$N$, full decorrelation, viz a decrease of the 
cross-correlation to essentially zero, occurs 
within the timescale of the initial exponential 
decorrelation, which is the hallmark of classical 
turbulent chaos. Decorrelation is a bit slower
for $N = 300$, which indicates that the kernel
series framework is close to a partially predictable
chaotic state in this regime.

\section{Generic delay distributions}\label{sect:generic}

So far we developed the theory for kernel
series frameworks that are associated with 
a delay system containing a single delay $T$.
Here we generalized our framework to systems 
characterized from the start by a 
pre-determined distribution of time delays.

As a first step we consider systems containing
a finite number $M$ of discrete time delays
$T_j$. Generically, distinct time delays 
could serve specific functional roles, like 
in a logistic equation with mixed
delays, $\dot{x}(t) = x(t-T_1)[1-x(t-T_2)]$.
Alternatively, as considered here, the 
history $h(t)$ entering (\ref{eq:general_dde}) 
is assumed to contain $M$ terms with
relative weights $\kappa_j\ge0$,
\begin{equation}
h(t) = \sum_j \kappa_j x(t-T_j),
\qquad\quad
\sum_j \kappa_j=1\,.
\label{eq_history_M_discrete}
\end{equation}
For every time delay $T_j$ one constructs
a kernel series of length $N_j$
\begin{equation}
\left\{ K_i^{(N_j, T_j)}(t) \right\}_{i = 1, \dots, N_j}, 
\qquad\quad j = 1, \dots, M\,,
\label{eq_kernel_M}
\end{equation}
as defined by (\ref{eq:kernel_gamma}).
In addition, one introduces $\sum_j N_j$
corresponding auxiliary variables obeying 
suitably generalized versions of (\ref{eq_dot_x_m}).

For a general distribution of time delays,
the history is given by
\begin{equation}
h(t) = \int_0^\infty  \kappa(\tau) x(t-\tau)\,\mathrm{d}\tau,
\qquad\quad
\int_0^\infty  \kappa(\tau)\,\mathrm{d}\tau = 1\,,
\label{eq_history_M_continuous}
\end{equation}
which generalizes (\ref{eq_history_M_discrete}).
Next one uses the fact that exponential shapelets
$\{ \psi_i(t; \beta) \}_{i \in \mathbb{N}}$
constitute an orthogonal basis on the positive
real axis \cite{berge2019exponential}. The
delay kernel can hence be expanded in shapelets,
as
\begin{equation}
\kappa(t) = \sum_i c_{i} \psi_i(t; \beta), 
\qquad\quad 
c_i = \int_0^\infty \psi_i(t; \beta) \kappa(t)\,\mathrm{d}t\,,
\label{eq_kappa_expansion}
\end{equation}
where the scale parameter $\beta$ can be used for 
optimization purposes, e.g.\ to minimize the number 
of non-zero expansion coefficients $c_j$.
Given that exponential shapelets are superpositions
of the elementary kernels introduced in (\ref{eq:kernel_def}),
one can convert (\ref{eq_kappa_expansion}) into
a kernel series framework, albeit at the cost of 
a possibly diverging number of auxiliary variables.

\section{Chaotic diffusion}\label{sec:chaotic_diffusion}
%
\begin{figure}[t]
  \centering
  \includegraphics[]{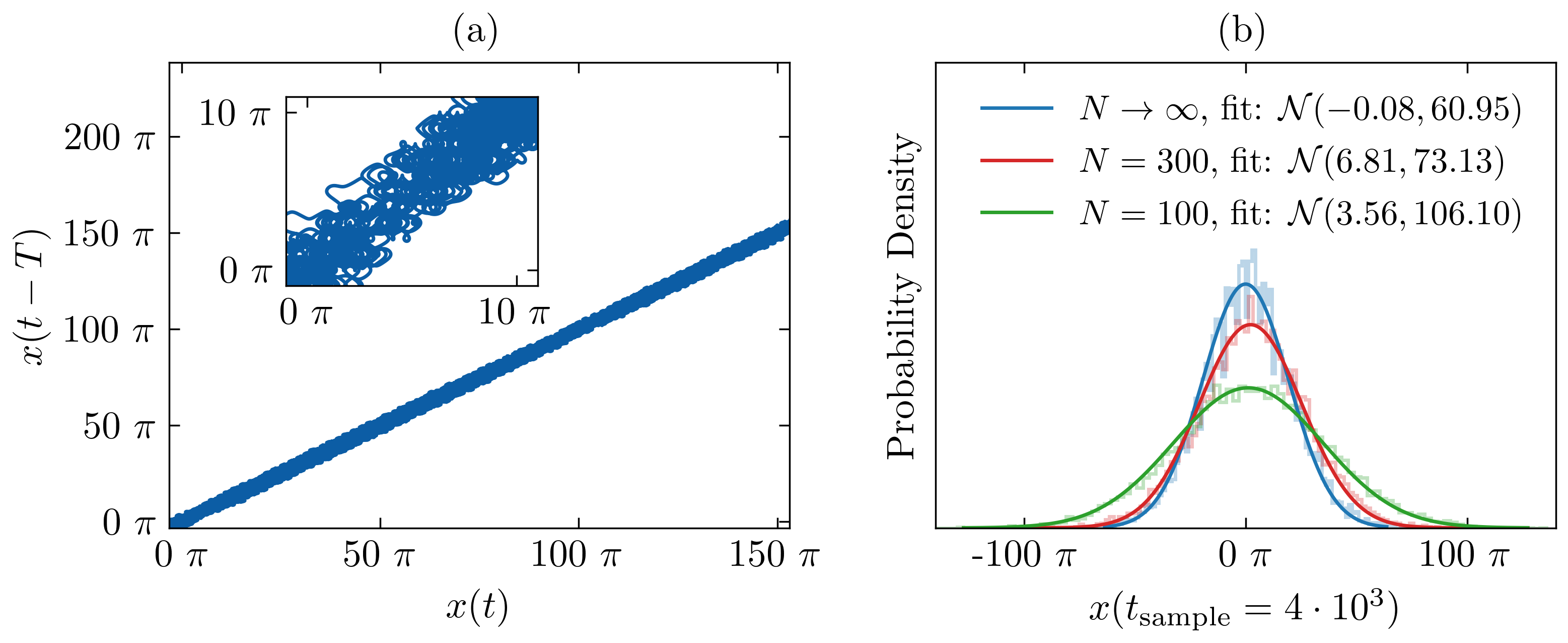}
  \caption{\textbf{Chaotic diffusion.} In \textbf{(a)}, for a time delay
  $T=20$, the stroboscopic projection of
  the chaotic attractor of (\ref{eq:ikeda}).
  The inset shows how the trajectory moves
  along the band of previously disconnected
  limit cycles.
  The trajectory spans the entire
  real axis, $-\infty < x(t) < \infty$.
  \textbf{(b)} The probability density 
  for the distance from the origin, also for
  $T=20$, when sampling at
  $t_\mathrm{sample} = 4\cdot 10^3$.
  The three data sets shown correspond to a
  direct integration of (\ref{eq:ikeda}) 
  (blue curve) and to the behavior obtained when
  using the kernel series framework, with 
  orders $N = 100, 300$ (green and red curve,
  respectively). A total of $10^4$ initial 
  conditions close to the origin have been
  considered. The data is fitted to a Gaussian 
  $\mathcal{N}(\mu, \sigma)$ with mean $\mu$ 
  and standard deviation $\sigma$, as given in the 
  legend. Chaotic diffusion is 
  more prominent for smaller orders $N$.}
  \label{fig:diffusive_chaos}
\end{figure}
So far we investigated the Mackey-Glass system 
within the kernel series framework. 
Next we turn towards a time-delayed 
system exhibiting chaotic diffusion. 
The flow,
\begin{equation}\label{eq:ikeda}
\dot{x}(t) = \sin\left(x(t-T)\right)\,,
\end{equation}
is given by a time-delayed sinusoidal function.
Possible additional parameters, like an overall 
prefactor, can be eliminated by rescaling $x$ and $t$.
The properties of (\ref{eq:ikeda}) have been studied 
in detail in \cite{sprott2007simple}.
The fixpoints $x^* = k\pi$, with $k\in\mathbb{N}$,
loose stability in a Hopf bifurcation which gives 
rise to a stable limit cycle when increasing the 
time delay $T$. When further increasing the time 
delay, the period of the limit cycle doubles, with
the system transitioning afterwards to a chaotic 
regime via an attractor-merging crisis. These 
features are reproduced when using the kernel 
series framework, with the accuracy increasing
as a function of $N$.

Linearizing (\ref{eq:ikeda}) yields 
(\ref{eq:linear_dde}) with $a = 1$ and 
$b = c =0$. Therefore, the discussion is 
qualitatively the same as for the Mackey-Glass
system, as given in section \ref{sec:local_regime}. 
Again, the relative differences between the loci 
of the bifurcation points in the time delayed
system compared to its representation in the 
kernel series framework scales as $1/N$, as
we did demonstrate previously for the 
Mackey-Glass system.

In \cite{sprott2007simple} it is noted that  (\ref{eq:ikeda}) shows chaotic diffusion, which
implies that the statistics of the
trajectory is given in the chaotic regime 
by 1D Brownian diffusion. This means that 
the position $x(t_\mathrm{sample})$ is 
normal-distributed when evolving the system 
up to a sampling time $t_\mathrm{sample} \gg 1$ 
for a number of random initial conditions
near the origin. One finds $\mu\to0$
for the mean \cite{sprott2007simple}.
Typical for diffusive behavior is a
linearly increasing variance,
$\sigma^2\sim t_\mathrm{sample}$, a
behavior that is possible because the
chaotic state of (\ref{eq:ikeda}) forms
a band of trajectories from the previously 
disconnected chain of limit cycles obtained
via $2\pi$-shifts \cite{sprott2007simple}.
As an illustration, the chaotic attractor of (\ref{eq:ikeda}) is
presented in figure~\ref{fig:diffusive_chaos}~(a)
for a time delay of $T = 20$. The trajectory 
spans the entire real axis, $-\infty < x(t) < \infty$ \cite{sprott2007simple}.

In figure~\ref{fig:diffusive_chaos}~(b) we 
numerically evaluate the distribution 
of $x(t_\mathrm{sample})$ for 
different orders $N = 100, 300$ of 
the kernel series framework and compare
to the case where (\ref{eq:ikeda}) is 
integrated directly ($N\to\infty$), using 
a sampling time of $t_\mathrm{sample} = 4\cdot 10^3$,
$10^4$ initial conditions and considering a time delay $T = 20$.
We find that
for different orders $N$ the data is fitted 
accurately by a Gaussian distribution. The 
mean remains close to zero for all cases considered. 
The standard deviation, which is directly 
related to the diffusion constant, is however
strongly dependent on the order $N$, growing
in size when lowering $N$. It would be interesting 
to investigate the scaling of the standard deviation 
with $1/N$, viz as a function of the variance of the
distribution of time delays, which is however
computationally demanding. We leave this aspect 
for future work. In any case our results show
that chaotic diffusion is substantially more
pronounced when the distribution of time delays
has a positive width.

\section{Discussion}\label{sec:discussion}

There are several reasons why distributed
time delays are important: Firstly, on a fundamental
level, because memory formation takes time,
as pointed out in the introduction. Secondly,
because natural processes, like the dynamics
of blood cells described by the
Mackey-Glass system \cite{mackey1977oscillation},
are often intrinsically noisy. For memories
of aggregate quantities, like the concentration
of cells, biological variability translates into
a corresponding distribution function. A similar
argument can be made for socio-economic and
climate models with delays 
\cite{keane2017climate,sportelli2022goodwin}, 
for which
delayed feedback is transmitted in general via
a cascade of intermediate processes.

Here we systematically investigated
the effects broadened memory kernels have on
the dynamics of time delayed systems by replacing
discrete delays with distributed delays.
It is known (e.g.\ \cite{hurtado2019generalizations,
cassidy2021distributed})
that a specific type of delay distributions,
gamma distributions with integer exponents, 
can be mapped exactly onto a set of $N+1$ 
ordinary differential equations. We here denoted this procedure
the `kernel series framework'.
Gamma distributions take the form of broadened
$\delta$-functions, as illustrated in 
figure\,\ref{fig:kernel_convergence},
which allows to recover the case of a
single time delay $T$ as the limiting
case $N\to\infty$. The kernel series
framework is hence well suited for
the systematic study of the influence of
distributed time delays on the dynamical
phase diagram. Alternatively one may use
the kernel series framework as an
approximation to a given delay differential
equation. From a computational point of view
one has to weigh the perks of the kernel series
framework against the necessity to solve a much
higher dimensional system of differential equations.
The complexity of differential equation integration
usually scales linearly as $\mathcal{O}(n)$, where
$n$ denotes the size of the system
\cite{ansmann2018efficiently}. If the required
order in the kernel series framework is high,
integration may be computationally demanding.

In this paper, we studied numerically the influence of time delay distributions
for the Mackey-Glass system as well as a simple time delayed system with
sinusoidal nonlinearity. Good agreement is observed in the local regime for the
stability of fixed points, for which we prove analytically that corrections
scale as $1/N$. We also find that higher-order phenomena, such as period
doubling transitions, the occurrence of chaotic dynamics and chaotic diffusion,
are substantially more sensitive to the introduction of distributed delays. It
may hence be difficult to compare the predictions of dynamical systems with
precisely defined time delays with observations, in particular when distributed
time delays play an important role in the respective real-world applications.

\section{Appendix}\label{sec:appendix}

\subsection{Hopf bifurcation in perturbation theory}

An analytical estimate the Hopf bifurcation point 
occurring within a kernel series framework of order $N$,
as defined by (\ref{eq:kernel_approx}), may be attained 
through perturbation theory. In section\,\ref{sec:local_regime}
we showed that the characteristic polynomial 
$P_N(\lambda)$ approaches
(\ref{eq:characteristic_eq}) as $(1+x/k)^k$ for $k\to\infty$. 
The relation
\begin{equation}
\label{eq:exp_restglied}
\left( 1 + \frac{x}{k} \right)^k = 
\mathrm{e}^x \left( 1 - \frac{x^2}{2k} \right)
\end{equation}
holds in leading order $1/k$, which is seen by taking 
the derivative of both sides and comparing
in leading order. Inserting (\ref{eq:exp_restglied}) 
into the characteristic polynomial
(\ref{eq:ode_characteristic_eq}), we obtain
\[
P_N(\lambda) = \lambda + b + a \mathrm{e}^{-\lambda T} - 
\epsilon \left[(b + \lambda) \frac{T^2 \lambda^2}{2}\right] =: 
P_\infty(\lambda) + C(\lambda, \epsilon) = 0\,,
\]
where $\epsilon = 1/N$ and $C(\lambda, \epsilon)$ 
denotes a perturbation of (\ref{eq:characteristic_eq}),
the characteristic equation $P_\infty(\lambda)$ of the
system containing only a single delay.
At the Hopf bifurcation point $T_\mathrm{Hopf}$,
the real part of 
$\lambda$ vanishes. Thus, we make the following
ansatz at the Hopf bifurcation point for 
$q := \mathrm{Im}(\lambda)$ and $T_\mathrm{Hopf}$ 
in terms of $\epsilon$:
\begin{equation}\label{eq:perturb_ansatz}
q = \sum_{i=0}^\infty c_i \epsilon^i, \qquad T_\mathrm{Hopf} = 
\sum_{i=0}^\infty k_i \epsilon^i.
\end{equation}
In zeroth order the solution (\ref{eq:hopf_T_c}) found for the
DDE is reproduced. In first order
we attain
\[
k_1 = \frac{\arccos^2(-b/a) \left(b \sqrt{a^2-b^2} + 
a^2 \arccos(-b/a)\right)}{2 (a^2 - b^2)^{3/2}}\,.
\]
For our usual values $a=0.4$ and $b=0.1$ we thus 
find in first order perturbation theory
\begin{equation}
T_\mathrm{Hopf} \simeq  4.708 + 
9.458\cdot\epsilon,
\qquad\quad
\epsilon = \frac{1}{N}\,.
\label{eq:appedix_T_Hopf}
\end{equation}
j
\section*{Data availability statement}

All data that support the findings of this 
study are included within the article 
(and any supplementary files).

\section*{Acknowledgments}

This project has received funding from the European Union's Horizon 2020
research and innovation programme under grant agreement No 101016233.

\noindent The authors thank Michael C.~Mackey, Bulcs\'u S\'andor and Rahil Valani
for inspiring comments.

\section*{ORCID ID}

C. Gros \href{ https://orcid.org/0000-0002-2126-0843}{https://orcid.org/0000-0002-2126-0843}\\
\noindent D. H. Nevermann \href{https://orcid.org/0000-0002-4607-5142} {https://orcid.org/0000-0002-4607-5142}

\section*{References}  
\bibliographystyle{unsrt}
\bibliography{paper.bib}

\end{document}